\documentclass[12pt]{amsart}
\usepackage[a4paper,hmargin=3cm,vmargin=3cm]{geometry}
\usepackage{amsfonts,amssymb,amscd,amstext,amsthm}
\usepackage{graphicx}
\usepackage[T1]{fontenc}
\usepackage[english]{babel}
\usepackage[utf8]{inputenc}
\usepackage[bookmarksnumbered,plainpages]{hyperref}
\usepackage{indentfirst}
\usepackage{color}
\usepackage{url}
\usepackage{comment}

\usepackage[shortlabels]{enumitem}
\usepackage{imakeidx}
\usepackage{mathrsfs}
\usepackage{stmaryrd}
\usepackage{textfit}
\usepackage{lipsum}
\usepackage{mathtools}
\usepackage{times}

\newcommand{\Ric}{\mbox{Ric}}

\newcommand{\tr}{\mbox{tr}}
\newcommand{\spa}{\mbox{span}}
\newcommand{\Hess}{\mbox{Hess}}
\newcommand{\End}{\mbox{End}}
\newcommand{\grad}{\mbox{grad}}
\newcommand{\Id}{\mbox{Id}}

\newtheorem{theorem}{Theorem}[section]
\newtheorem{lemma}[theorem]{Lemma}
\theoremstyle{definition}

\newtheorem{example}[theorem]{Example}

\newtheorem{proposition}[theorem]{Proposition}

\theoremstyle{remark}

\numberwithin{equation}{section}

\newcommand{\Addresses}{{
\bigskip 
\footnotesize
Institute of Mathematics and Computer Science \newline
University of S\~ao Paulo, S\~ao Carlos, Brazil \\
\textit{E-mail address:} \texttt{mateusrodrigues@alumni.usp.br}
}}

\title{On submanifolds with constant M\"{o}bius Ricci curvature}

\date{}

\author{M. S. R. Antas}

\begin{document}

\begin{abstract}
In the M\"{o}bius geometry of submanifolds there exists an important M\"{o}bius invariant tensor called the M\"{o}bius metric which together with the M\"{o}bius second fundamental form, the Blaschke tensor and the M\"{o}bius form constitute a complete system of conformal invariants for umbilic-free submanifolds when the dimension $n \geq 3$. An umbilic-free immersion $f:M^{n}\to \mathbb{R}^{m}$ is said to have \emph{constant M\"{o}bius Ricci curvature} if its M\"{o}bius metric has constant Ricci curvature. In this paper, we classify umbilic-free immersions $f:M^{n}\to \mathbb{R}^{m}$ with constant M\"{o}bius Ricci curvature, flat normal bundle and two distinct principal normal vector fields. We prove that a submanifold with constant M\"{o}bius Ricci curvature and flat normal bundle has closed M\"{o}bius form, and then we classify such submanifolds provided the M\"{o}bius form vanishes identically on $M^{n}$. 
\end{abstract}

\maketitle

{\small
\noindent {\textbf{Mathematics Subject Classification:}}\hspace*{0.1cm}
Primary 53B25, 53C30, Secondary 53C40.
}
\vspace{.2cm} \\
{\small
\noindent {\textbf{Keywords:}}\hspace*{0.1cm}
M\"{o}bius metric, Einstein submanifolds, flat normal bundle, M\"{o}bius form, principal normals.
}

\section{Introduction}

The classical paper \cite{Wang} introduces two important M\"{o}bius invariant tensors named the \emph{M\"{o}bius metric} $\langle \cdot, \cdot \rangle^{*}$ and the \emph{M\"{o}bius second fundamental form} $\beta$ associated to an umbilic-free immersion $f:M^{n}\to \mathbb{R}^{m}$ (see \eqref{M\"{o}bius-metric} and \eqref{M\"{o}bius-sec-form}). It is also proved in Wang's work that $f$ is solely determined by $\langle \cdot, \cdot \rangle^{*}, \beta$ and the normal connection $\nabla^{\perp}$, up to M\"{o}bius transformations in $\mathbb{R}^{m}$, when $n \geq 3$ (see also Theorem 9.22 in \cite{DT}). The corresponding conformal Gauss, Codazzi and Ricci equations involve two other important M\"{o}bius invariant tensors named the \emph{Blaschke tensor} and the \emph{M\"{o}bius form} (see \eqref{blaschke} and \eqref{M\"{o}bius-form}).

Wang's work motivated several authors to study umbilic-free immersions whose M\"{o}bius metric has constant curvatures (see, for example, \cite{AT, sectional, Ricci, scalar, deformations, willmore, notes}).

An important class of umbilic-free submanifolds $f:M^{n}\to \mathbb{R}^{m}$ are those with \emph{constant M\"{o}bius curvature}, which means that $(M^{n}, \langle \cdot, \cdot \rangle^{*})$ has constant sectional curvatures. In \cite{sectional, deformations} was obtained the classification of umbilic-free hypersurfaces $f:M^{n}\to \mathbb{R}^{n+1}$, $n \geq 4$, with constant M\"{o}bius curvature. The proof of this theorem still holds for $n=3$, provided that the hypersurface has two distinct principal curvatures. In \cite{willmore} was obtained the classification of hypersurfaces $f:M^{3}\to \mathbb{R}^{4}$ that has constant M\"{o}bius curvature and three distinct principal curvatures. In \cite{AT}, the authors classify the submanifolds $f:M^{n}\to \mathbb{R}^{n+p}$, $n \geq 5$ and $2p \le n$ with constant M\"{o}bius curvature and flat normal bundle. The proof of this theorem also works for submanifolds of dimension $n=4$ that have at least one principal normal vector field of multiplicity greater than one. It is also proved in \cite{AT} a classification of the submanifolds with constant M\"{o}bius curvature, flat normal bundle and two distinct principal normal vector fields with no restriction on the codimension.

Another notable class of umbilic-free submanifolds $f:M^{n}\to \mathbb{R}^{m}$ are those with \emph{constant M\"{o}bius Ricci curvature}, meaning that $(M^{n}, \langle \cdot, \cdot \rangle^{*})$ has constant Ricci curvatures. It is obvious that all submanifolds with constant M\"{o}bius sectional curvature have constant M\"{o}bius Ricci curvature, but the converse might be false if $n \geq 4$. However, for $n=3$, submanifolds with constant M\"{o}bius Ricci curvature have constant M\"{o}bius sectional curvature. In the paper \cite{Ricci}, the authors classify the hypersurfaces $f:M^{n}\to \mathbb{R}^{n+1}$ with constant M\"{o}bius Ricci curvature and two distinct principal curvatures. Precisely, they showed that for $n \geq 4$, $f(M^{n})$ is locally M\"{o}bius equivalent to a hypersurface with constant M\"{o}bius sectional curvature or the inverse image under the stereographic projection of the torus $\mathbb{S}^{k}\left(\sqrt{\frac{k-1}{n-2}}\right)\times \mathbb{S}^{n-k}\left(\sqrt{\frac{n-k-1}{n-2}}\right)\subset \mathbb{S}^{n+1},$ $1 <k<n-1.$ Very recently, Li and Xie proved in \cite{notes} that there are not hypersurfaces of dimension $n \geq 4$ with constant M\"{o}bius Ricci curvature and at least three distinct principal curvatures.

In this paper, we initiate the study of umbilic-free submanifolds $f:M^{n}\to \mathbb{R}^{m}$, $n \geq 3$, with constant M\"{o}bius Ricci curvature and higher codimension. We restrict ourselves to the case in which the submanifold has flat normal bundle. Our first main result is a classification of these submanifolds with two distinct principal normal vector fields. 
\begin{theorem}\label{main1}
    Let $f:M^{n}\to \mathbb{R}^{n+p}$, $n \geq 3$, be an isometric immersion with flat normal bundle and constant M\"{o}bius Ricci curvature that has two distinct principal normal vector fields. Then $f(M^{n})$ is, up to a M\"{o}bius transformation in $\mathbb{R}^{n+p},$ an open subset of one of the following submanifolds
      \begin{description}
        \item[a)] a submanifold with constant M\"{o}bius sectional curvature and flat normal bundle (see Theorem 4.1 of \cite{AT});
        \item[b)] the inverse image under the stereographic projection of the torus $$\mathbb{S}^{k}\left(\sqrt{\frac{k-1}{n-2}}\right)\times \mathbb{S}^{n-k}\left(\sqrt{\frac{n-k-1}{n-2}}\right)\subset \mathbb{S}^{n+1},\quad 1 <k<n-1.$$ 
    \end{description}
\end{theorem}

In Lemma \ref{M\"{o}bius-closed} is proved that an umbilic-free submanifold $f:M^{n}\to \mathbb{R}^{m}$ with constant M\"{o}bius Ricci curvature and flat normal bundle has closed M\"{o}bius form. Our second main result consists in the following theorem.

\begin{theorem}\label{main22}
    Let $f:M^{n}\to \mathbb{R}^{m}$, $n \geq 3,$ be an umbilic-free isometric immersion with constant M\"{o}bius Ricci curvature, flat normal bundle and vanishing M\"{o}bius form. Then $f(M^{n})$ is locally M\"{o}bius equivalent to one of the following submanifolds 

    \begin{description}
         \item[a)]  $\mathbb{S}^{1}(r_1)\times \ldots \times \mathbb{S}^{1}(r_{n-s})\times \mathbb{R}^{s}\subset \mathbb{R}^{2n-s},$ for all $1 \le s \le n-1$.
         \item[b)]  the three dimensional cone over the flat torus $\mathbb{S}^{1}(r)\times \mathbb{S}^{1}(\sqrt{1-r^{2}})\subset \mathbb{S}^{3}.$
        \item[c)] the inverse image under the stereographic projection of the torus $$\mathbb{S}^{k}\left(\sqrt{\frac{k-1}{n-2}}\right)\times \mathbb{S}^{n-k}\left(\sqrt{\frac{n-k-1}{n-2}}\right)\subset \mathbb{S}^{n+1},\quad 1 <k<n-1.$$  
    \end{description}
\end{theorem}

Classifying a certain class of umbilic-free submanifolds with vanishing M\"{o}bius form is a topic that has attracted much attention. In \cite{scalar}, the authors gave a rigidity of the submanifolds with constant M\"{o}bius scalar curvature under the condition vanishing M\"{o}bius form. In \cite{isotropic}, the authors use the fact that M\"{o}bius form vanishes identically on $M^{n}$ to give a unified characterization of the minimal submanifolds with constant scalar curvature in space forms. If the M\"{o}bius form vanishes identically in $M^{n}$, the previous result was generalized in \cite{mean-curvature-scalar}, precisely the authors gave a characterization of submanifolds with parallel mean curvature vector field and constant scalar curvature in space forms. In \cite{veronese}, it is given a M\"{o}bius characterization of Veronese surfaces in $\mathbb{S}^{n}$, assuming the vanishing of the M\"{o}bius form. In \cite{A}, the author studies conformally flat submanifolds with vanishing M\"{o}bius form.

Since there are no hypersurfaces of dimension $n \geq 4$ with constant M\"{o}bius Ricci curvature and at least three distinct principal curvatures, the only example of hypersurface with constant M\"{o}bius Ricci curvature that is not constant M\"{o}bius sectional curvature is item b) of Theorem \ref{main1}, which has vanishing M\"{o}bius form and two distinct principal curvatures. Then it is important to know if there exist other examples of umbilic-free submanifolds with constant M\"{o}bius Ricci curvature and flat normal bundle with vanishing M\"{o}bius form. Our answer given in Theorem \ref{main22} shows that there are no such submanifolds.

This paper is organized as follows. In section \ref{preliminaries}, we review the basic invariants of submanifolds in $\mathbb{R}^{m}$. In section \ref{two}, we prove Theorem \ref{main1} and in section \ref{three} we prove Theorem \ref{main22}.

\section{Preliminaries}\label{preliminaries}

Following the seminal paper \cite{Wang}, we recall the M\"{o}bius invariant tensors and the compatibility equations associated to an umbilic-free submanifold $f:M^{n}\to \mathbb{R}^{m}$.

The function $$\rho^{2}=\frac{n}{n-1}(||\alpha^{f}||^{2}-n||\mathcal{H}^{f}||^{2})$$ does not vanish on $M^{n}$, where $\mathcal{H}^{f}$ and $||\alpha^{f}||$ stand for the mean curvature vector field of $f$ and the norm of the second fundamental form of $f$, respectively. In the following, we define the M\"{o}bius lift of an umbilic-free immersion $f:M^{n}\to \mathbb{R}^{m}$. We denote by $\mathbb{L}^{m+2}$ the \emph{Minkowski} space of dimension $m+2$. The \emph{upper light cone} $\mathbb{V}^{m+1}_{+}$ of $\mathbb{L}^{m+2}$ is the upper connected component of the subset $$\mathbb{V}^{m+1}=\{v \in \mathbb{L}^{m+2}: \langle v,v \rangle=0\,\, \mbox{and}\,\,\, v \ne 0\}.$$ Nevertheless, when restricted to its subset $$\mathbb{E}^{m}=\mathbb{E}^{m}_{w}=\{p \in \mathbb{V}^{m+1}_{+}: \langle p, w \rangle=1\},$$ it gives rise to a model of the $m$-dimensional Euclidean space $\mathbb{R}^{m}$ for each fixed $w \in \mathbb{V}^{m+1}_{+}$. Indeed, for any $p_0\in \mathbb{E}^{m}$ and any linear isometry $C:\mathbb{R}^{m}\to (\spa\{p_0,w\})^{\perp}\subset \mathbb{L}^{m+2}$, the map $\Psi=\Psi_{p_0,w,C}: \mathbb{R}^{m}\to \mathbb{L}^{m+2}$ given by $$\Psi(x)=p_0+Cx-\frac{1}{2}||x||^{2}w,$$ is an isometric embedding such that $\Psi(\mathbb{R}^{m})=\mathbb{E}^{m}.$ The \emph{M\"{o}bius lift} of $f$ is the immersion $F:M\to \mathbb{V}^{m+1}_{+}\subset \mathbb{L}^{m+2}$ given by $F=\rho(\Psi\circ f).$ Since its differential is $F_{*}X=X(\rho)\Psi \circ f+\rho \Psi_{*}\circ f_{*}X,$ it follows from the fact that $\Psi$ is a light-like normal vector field that the induced metric by $F$ is \begin{equation}\label{M\"{o}bius-metric}\langle \,\,,\,\rangle^{*}=\rho^{2}\langle \,\,,\,\rangle_{f},\end{equation} which is called the \emph{M\"{o}bius metric} determined by $f$, where $\langle \, \, , \,\rangle_{f}$ denotes the metric induced by $f$. The second fundamental form of $F$ is given by \begin{equation}\label{sf}\alpha^{F}(X,Y)=\phi(\beta(X,Y))-\psi(X,Y)F-\langle X,Y\rangle^{*}\zeta\end{equation} where $\beta$ and $\psi$ are tensors defined below, $\phi$ is a vector bundle isometry of the normal bundle $N_fM$ onto a subbundle $V$ of $N_FM$ and $\zeta$ is a vector field on $M^{n}$ such that $\{F, \zeta\}$ forms a pseudo-orthonormal frame in $V^{\perp}$ with $$\langle F,F\rangle=0=\langle \zeta, \zeta \rangle \quad \mbox{and} \quad \langle F, \zeta \rangle=1.$$ (See Proposition 9.17 of \cite{DT} for more details.)

\emph{The M\"{o}bius second fundamental form} of $f:M^{n}\to \mathbb{R}^{m}$ is the symmetric section $\beta=\beta^{f} \in \mbox{Hom}^{2}(TM, TM; N_{f}M)$ defined by \begin{equation} \label{M\"{o}bius-sec-form}\beta(X,Y)=\rho (\alpha^{f}(X,Y)-\langle X,Y\rangle \mathcal{H}^{f})\end{equation} 
 for all  $ X,Y \in \mathfrak{X}(M)$. Notice that $\beta$ is traceless and that its norm  with respect to  $\langle \,,\,\rangle^{*}$ is $||\beta||_{*}=\sqrt{\frac{n-1}{n}}$. 

The \emph{Blaschke tensor} $\psi=\psi^{f}$ of $f$ is the symmetric $C^{\infty}(M)$-bilinear form given by \begin{equation}\label{blaschke}\psi(X,Y)=\frac{1}{\rho}\langle \beta(X,Y), \mathcal{H}\rangle+\frac{1}{2\rho^{2}}(||\grad^{*}\rho||^{2}_{*}+||\mathcal{H}||^{2})\langle X,Y\rangle^{*}-\frac{1}{\rho}\Hess^{*}\rho(X,Y)\end{equation} and its \emph{M\"{o}bius form} $\omega=\omega^{f}$ is the normal bundle valued one-form defined by \begin{equation}\label{M\"{o}bius-form}\omega(X)=-\frac{1}{\rho}(\nabla^{\perp}_{X}\mathcal{H}^{f}+\beta(X,\grad^{*}\rho))\end{equation} where $\grad^{*}$ and $\Hess^{*}$ denote the gradient and Hessian on $(M^{n}, \langle \,\,,\,\rangle^{*})$.

 Associated with $\beta$ one defines the  \emph{M\"{o}bius third fundamental form } $III_{\beta} : \mathfrak{X}(M)\times \mathfrak{X}(M)\to \mathbb{R}$  by \begin{equation*}\label{thirdform}
	III_{\beta}(X,Y)=\sum_{i=1}^{n}\langle \beta(X,X_i), \beta(Y,X_i)\rangle,
\end{equation*} where  $X_1, \ldots, X_n$ is an orthonormal frame with respect to $\langle \, \,,\,\rangle^{*}.$

\begin{proposition}[\cite{Wang}]\label{blaschke2}
\textit{The Blaschke tensor is given in terms of the M\"{o}bius metric and the M\"{o}bius third fundamental form by }  \begin{equation*}
(n-2)\psi(X,Y)=\Ric^{*}(X,Y)+III_{\beta}(X,Y)-\frac{n^{2}s^{*}+1}{2n}\langle X,Y\rangle^{*},
\end{equation*} \textit{for all $X,Y \in \mathfrak{X}(M)$, where $\Ric^{*}$ and $s^{*}=\frac{1}{n(n-1)}\tr\, \Ric^{*}$ are the Ricci curvature and the scalar curvature of $(M^{n}, \langle \, \,,\,\rangle^{*}).$ In particular}, \begin{equation*}\tr\, \psi=\frac{n^{2}s^{*}+1}{2n}. \end{equation*}
\end{proposition}

The M\"{o}bius metric, the M\"{o}bius second fundamental form, the Blaschke tensor and the M\"{o}bius form of $f$ are M\"{o}bius invariant tensors that satisfy the following compatibility equations.

\begin{proposition}[\cite{Wang}]\label{prop1}
	\textit{The following equations hold:}
	
	\begin{description}
		\item[(i)] \textit{The conformal Gauss equation}: \begin{align}\label{G}
			\langle R^{*}(X,Y)Z,W\rangle^{*}=&\langle \beta(X,W), \beta(Y,Z)\rangle - \langle \beta(X,Z), \beta(Y,W)\rangle\nonumber\\ 
			&+\psi(X,W)\langle Y,Z \rangle^{*}+\psi(Y,Z)\langle X,W\rangle^{*}\nonumber\\
			&-\psi(X,Z)\langle Y,W \rangle^{*}-\psi(Y,W)\langle X,Z \rangle^{*}
		\end{align} \textit{for all} $X,Y,Z,W \in \mathfrak{X}(M)$.
		\item[(ii)] \textit{The conformal Codazzi equations}:
		\begin{equation}\label{C1}
			({^{f}\nabla^{\perp}_{X}}\beta)(Y,Z)-({^{f}\nabla^{\perp}_{Y}}\beta)(X,Z)=\omega((X \wedge Y)Z)
		\end{equation} \textit{and} \begin{equation}\label{C2}
			(\nabla^{*}_{X}\psi)(Y,Z)-(\nabla^{*}_{Y}\psi)(X,Z)=\langle \omega(Y), \beta(X,Z)\rangle-\langle \omega(X), \beta(Y,Z)\rangle
		\end{equation} \textit{for all} $ X,Y,Z \in \mathfrak{X}(M),$ \textit{where} $$ ({^{f}\nabla^{\perp}_{X}}\beta)(Y,Z)={^{f}\nabla^{\perp}_{X}}\beta(Y,Z)-\beta(\nabla^{*}_{X}Y,Z)-\beta(Y, \nabla^{*}_{X}Z),$$  $$(\nabla^{*}_{X}\psi)(Y,Z)=X(\psi(Y,Z))-\psi(\nabla^{*}_{X}Y, Z)-\psi(Y, \nabla^{*}_{X}Z)$$ \textit{and}
		$$(X \wedge Y)Z=\langle Y,Z\rangle^{*} X-\langle X,Z\rangle^{*} Y.$$ 
		\item[(iii)] \textit{The conformal Ricci equations}: \begin{equation}\label{R1}
			d\omega(X,Y)=\beta(Y, \hat{\psi}X)-\beta(X, \hat{\psi}Y)
		\end{equation} \textit{and} \begin{equation}\label{R2}
			\langle R^{\perp}(X,Y)\xi, \eta \rangle=\langle [B_{\xi}, B_{\eta}]X,Y\rangle^{*}
		\end{equation} \textit{for all} $ X,Y \in \mathfrak{X}(M)$ \textit{and} $\xi, \eta \in \Gamma(N_{f}M)$, \textit{with} $\hat{\psi}, B_{\xi} \in \Gamma(\End(TM))$ \textit{given by} \begin{equation*}
			\langle \hat{\psi}X,Y\rangle^{*}=\psi(X,Y) \quad \mbox{and} \quad \langle B_{\xi}X,Y\rangle^{*}=\langle \beta(X,Y), \xi \rangle.
		\end{equation*}
	\end{description}	
\end{proposition}

The next result is a fundamental theorem.
\begin{theorem}[\cite{Wang}]
    \textit{Let $f,g:M^{n}\to \mathbb{R}^{m}$, $n \geq 3$, be immersions free of umbilical points. Then there exists a conformal transformation $\tau:\mathbb{R}^{m}\to \mathbb{R}^{m}$ such that $g=\tau \circ f$ if and only if $f$ and $g$ share the same M\"{o}bius metric and there exists a vector bundle isometry $\mathcal{T}:N_{f}M \to N_{g}M$ such that $$\mathcal{T} {^{f}\nabla^{\perp}}={^{g}\nabla^{\perp}}\mathcal{T} \quad \mbox{and} \quad  \mathcal{T} \circ \beta^{f}=\beta^{g}.$$}
\end{theorem}

Let $f:M^{n}\to \mathbb{R}^{n+p}$ be a proper isometric immersion with flat normal bundle and let ${\eta}_1, \ldots, {\eta}_k$, $k \geq 2$, be the distinct principal normal vector fields of $f$ whose associated eigenbundles are denoted by $E_{\eta_1}, \ldots, E_{\eta_k}.$ Given $X\in \Gamma(E_{\eta_i})$ and $Y \in \mathfrak{X}(M)$, we have \begin{align}
    \beta(X,Y)&=\rho(\alpha(X,Y)-\langle X,Y\rangle \mathcal{H})\nonumber\\
    &=\rho(\langle X,Y\rangle \eta_i-\langle X,Y\rangle \mathcal{H})\nonumber\\
    &=\langle X,Y\rangle^{*}\rho^{-1}(\eta_i-\mathcal{H}).\label{M\"{o}bius-principal}
\end{align}

We call the normal vector fields $\bar{\eta}_i=:\rho^{-1}(\eta_i-\mathcal{H})$, $1 \le i \le k,$ by \emph{M\"{o}bius principal normal vector fields} of $f.$ Notice that $E_{\eta_i}$ is the smooth distribution associated with $\bar{\eta}_i$ and the number of distinct M\"{o}bius principal normal vector fields is $k$.

An umbilic-free isometric immersion $f:M^{n}\to \mathbb{R}^{m}$ is said to have constant M\"{o}bius Ricci curvature $\lambda$ if $\Ric^{*}(X,Y)=\lambda \langle X,Y\rangle^{*}$ for all $X,Y \in \mathfrak{X}(M),$ where $\Ric^{*}$ denotes the ricci tensor of the M\"{o}bius metric. We point out that this condition is a natural generalization of the constant sectional curvature of $(M^{n}, \langle\cdot, \cdot \rangle^{*})$. Moreover, $f$ has constant M\"{o}bius Ricci curvature $\lambda$ if and only if its M\"{o}bius lift $F$ has Ricci curvature $\lambda.$  It is also well known that all three dimensional Einstein manifold is of constant sectional curvature.

\section{The case of two distinct principal normals}\label{two}

An interesting property of an umbilic-free submanifold $f:M^{n}\to \mathbb{R}^{n+p}$, $n \geq 3$, with constant M\"{o}bius Ricci curvature and flat normal bundle is that it has closed M\"{o}bius form. In particular, there exists an orthonormal frame $\{X_1, \ldots, X_n\}$ that diagonalizes $\beta$ and $\psi$ 
simultaneously. We fix this frame through this section.

 \begin{lemma}\label{M\"{o}bius-closed}
		Let $f:M^{n}\rightarrow \mathbb{R}^{n+p}$ be an isometric immersion with constant M\"{o}bius Ricci curvature $\lambda$ and flat normal bundle. Then the M\"{o}bius form of $f$ is closed.
	\end{lemma}
 \begin{proof}
     The conformal Ricci equation \eqref{R2} shows that there exists an orthonormal frame $\mathcal{B}:=\{X_1, \ldots, X_n\}$ with respect to the M\"{o}bius metric such that $\beta(X_i,X_j)=0,$ $i \ne j.$
  By hypothesis, $\Ric^{*}(X,Y)=\lambda\langle X,Y\rangle^{*}$ for all $X,Y \in \mathfrak{X}(M).$ Hence $\mathcal{B}$ diagonalizes $\Ric^{*}$. 
     Using the relation in Proposition \ref{blaschke2}, we conclude that the Blaschke tensor $\psi$ is also diagonalizable by $\mathcal{B}$. Therefore, $\omega$ is closed  by \eqref{R1}.
 \end{proof}

The goal of this section is to prove Theorem \ref{main1}. In order to prove this result, we need the following result, which is proved in Lemma 4.1 of \cite{Antas1}.

\begin{lemma}\label{class-two}
	Let $f:M^{n}\to \mathbb{R}^{n+p}$, $n \geq 3$ and $p \geq 1$, be an isometric immersion with flat normal bundle that has two distinct principal normal vector fields $\eta_1$ and $\eta_2$ of multiplicities $k$ and $n-k$, for $1 \le k \le n-2$, respectively. If $\omega^{f}$ is closed then $f(M^{n})$ is M\"{o}bius equivalent to an open subset of a cylinder, a generalized cone or a rotational submanifold over an isometric immersion $g:M^{k}\to \mathbb{Q}^{p+k}_{c}$, where $c=0,1$ or $-1$, respectively. In particular, if $k \geq 2$ then $g$ is umbilical.
\end{lemma}

In the following, we describe the M\"{o}bius invariant tensors of the examples of the preceding lemma when $k=1.$

\begin{example}\label{cylinder}
	The cylinder in $\mathbb{R}^{n+p}$ over a curve $\gamma:I \to \mathbb{R}^{p+1}$ is defined by $$f=\gamma \times \Id:I \times \mathbb{R}^{n-1}\to \mathbb{R}^{n+p},$$ where $\Id:\mathbb{R}^{n-1}\to \mathbb{R}^{n-1}$ is the identity map.

    The M\"{o}bius metric determined by $f$ is given by $$\langle \,\, , \, \rangle^{*}=\kappa^{2}(ds^{2}+du_2^{2}+\ldots+du_n^{2}),$$ where $s$ is the arc-lenght parameter of $\gamma$, $\kappa(s)$ is the first curvature function of $\gamma$ and $(u_2, \ldots, u_n)$ are canonical coordinates on $\mathbb{R}^{n-1}.$

    The M\"{o}bius second fundamental form and the Blaschke tensor of $f$ are written in the orthonormal basis $\{\frac{1}{\kappa}\frac{d}{ds}, \frac{1}{\kappa}\frac{d}{du_2}, \ldots, \frac{1}{\kappa}\frac{d}{du_n}\}$ as \begin{align*}
		\beta\left(\kappa^{-1} \frac{d}{ds}, \kappa^{-1}\frac{d}{ds}\right)&=\frac{n-1}{n}\xi \quad \mbox{and} \quad
		\beta \left(\kappa^{-1}\frac{d}{du_i}, \kappa^{-1}\frac{d}{du_i}\right)=-\frac{1}{n}\xi,\\
        \psi\left(\kappa^{-1}\frac{d}{ds}, \kappa^{-1}\frac{d}{ds}\right)&=-\frac{\kappa_{ss}}{\kappa^{3}}+\frac{3}{2}\frac{k_s^{2}}{\kappa^{4}}+\frac{2n-1}{2n^{2}},\\ 
        \psi\left(\kappa^{-1}\frac{d}{du_i}, \kappa^{-1}\frac{d}{du_i}\right)&=-\frac{1}{2}\left(\frac{\kappa_s^{2}}{\kappa^{4}}+\frac{1}{n^{2}}\right),
	\end{align*} for all $2 \le i \le n$, where $\xi$ is a unit normal vector field of $f$.
\end{example}    

\begin{example}\label{cone1}
	The generalized cone in $\mathbb{R}^{n+p}$ over a curve $\gamma:I \to \mathbb{S}^{p+1}$ is the immersion $f:I \times \mathbb{H}^{n-1}\to  \mathbb{R}^{n+p}$ given by $$f(s, z)=\Theta \circ (\gamma, \Id)(s,z)=(z_1 \gamma(s), z_2,\ldots,z_{n-1}),$$ where $\mathbb{H}^{n-1}=\mathbb{R}_{+}\times \mathbb{R}^{n-2}$ is endowed with the hyperbolic metric $dz^{2}=\frac{1}{z_1^{2}}(dz_1^{2}+\ldots+dz_{n-1}^{2}),$ $Id:\mathbb{H}^{n-1}\to \mathbb{H}^{n-1}$ is the identity map and $\Theta: \mathbb{S}^{p+1}\times \mathbb{H}^{n-1}\to \mathbb{R}^{n+p}$ is the conformal diffeomorphism defined by $\Theta(y,z)=(z_1y, z_2, \ldots, z_{n-1}).$ 

    The M\"{o}bius metric determined by $f$ is given by $$\langle \,\,,\,\rangle^{*}=\kappa^{2}(ds^{2}+dz^{2}),$$ where $s$ is the arc-length parameter of $\gamma$ and $\kappa(s)$ is the first curvature function of $\gamma.$

    The M\"{o}bius second fundamental form and the Blaschke tensor of $f$ are written in the orthonormal basis $\{\kappa^{-1}\frac{d}{ds}, \kappa^{-1}z_1\frac{d}{dz_1}, \ldots, \kappa^{-1}z_1\frac{d}{dz_{n-1}}\}$ as 

    \begin{align*}
        \beta \left(\kappa^{-1}\frac{d}{ds}, \kappa^{-1}\frac{d}{ds}\right)&=\frac{n-1}{n}\xi \quad \mbox{and} \quad 
		\beta \left(\kappa^{-1}z_1\frac{d}{dz_i}, \kappa^{-1}z_1\frac{d}{dz_i}\right)=-\frac{1}{n}\xi,\\
        \psi \left(\kappa^{-1} \frac{d}{ds}, \kappa^{-1}\frac{d}{ds}\right)&=-\frac{\kappa_{ss}}{\kappa^{3}}+\frac{3}{2}\frac{\kappa_s^{2}}{\kappa^{4}}+\frac{1}{2\kappa^{2}}+\frac{2n-1}{2n^{2}},\\
        \psi \left(\kappa^{-1}z_1 \frac{d}{dz_i}, \kappa^{-1}z_1 \frac{d}{dz_i}\right)&=-\frac{1}{2}\left(\frac{\kappa_s^{2}}{\kappa^{4}}+\frac{1}{\kappa^{2}}+\frac{1}{n^{2}}\right),
    \end{align*} for all $1 \le i \le n-1$, where $\xi$ is a unit normal vector field of $f.$
\end{example}

\begin{example}\label{rotacao}
	The rotational submanifold over a curve $\gamma:I \to \mathbb{H}^{p+1}$ is the immersion $f:I \times \mathbb{S}^{n-1} \to \mathbb{R}^{n+p}$ defined by $$f(s,y)=\Theta \circ (\gamma, \Id)(s,y)=(z_1(s),\ldots, z_{p}(s), z_{p+1}(s)y),$$ where $\mathbb{H}^{p+1}=\{(z_1, \ldots, z_{p+1}): z_{p+1}>0\}$ is endowed with the hyperbolic metric $ds^{2}=\frac{1}{z_{p+1}^{2}}(dz_1^{2}+\ldots+dz_{p+1}^{2}),$ $\Id:\mathbb{S}^{n-1}\to \mathbb{S}^{n-1}$ is the identity map and $\Theta: \mathbb{H}^{p+1}\times \mathbb{S}^{n-1} \to \mathbb{R}^{n+p}$ is the conformal diffeomorphism $\Theta(z,y)=(z_1, \ldots, z_p, z_{p+1}y).$

    The M\"{o}bius metric determined by $f$ is $$\langle \,\,,\,\rangle^{*}=\kappa^{2}(\gamma^{*}ds^{2}+dy^{2}|_{\mathbb{S}^{n-1}}),$$ where $s$ is the arc-length parameter of $\gamma$ and $\kappa(s)$ is the first curvature function of $\gamma$.

    The M\"{o}bius second fundamental form and the Blaschke tensor of $f$ are written in the orthonormal basis $\{\kappa^{-1}\frac{d}{ds}, \kappa^{-1}e_1, \ldots, \kappa^{-1}e_{n-1}\}$ as \begin{align*}
	\beta \left(\kappa^{-1}\frac{d}{ds}, \kappa^{-1}\frac{d}{ds}\right)&=\frac{n-1}{n}\Theta_{*}\frac{\xi}{z_{p+1}} \quad \mbox{and} \quad 
	\beta \left(\kappa^{-1}e_i, \kappa^{-1}e_i\right)=-\frac{1}{n}\Theta_{*}\frac{\xi}{z_{p+1}},\\
    \psi \left(\kappa^{-1}\frac{d}{ds}, \kappa^{-1}\frac{d}{ds}\right)&=\frac{\kappa_{ss}}{\kappa^{3}}-\frac{5}{2}\frac{(\kappa_s)^{2}}{\kappa^{4}}-\frac{1}{2\kappa^{2}}+\frac{2n-1}{2n^{2}},\\
	\psi \left(\kappa^{-1}e_i, \kappa^{-1}e_i\right)&=-\frac{1}{2}\left[\frac{(\kappa_s)^{2}}{\kappa^{4}}-\frac{1}{\kappa^{2}}+\frac{1}{n^{2}}\right],
\end{align*} for all $1 \le i \le n-1$, where $e_1, \ldots, e_{n-1}$ is an orthonormal frame with respect to $dy^{2}|_{\mathbb{S}^{n-1}}$ and $\xi$ is a unit normal vector field of $\gamma.$
 \end{example}

Assume that the examples above have constant M\"{o}bius Ricci curvature $\lambda$ and denote by $\{X_1, \ldots,X_n\}$ their orthonormal frame with respect to $\langle \,\,,\,\rangle^{*}$. Since, by definition, $\Ric^{*}(X,Y)=\tr(Z \mapsto R^{*}(Z,X)Y)$, we have $$\lambda\langle X,Y\rangle^{*}=\sum_{i=1}^{n}\langle R^{*}(X_i,X)Y,X_i\rangle^{*}.$$ By the conformal Gauss equation, for each $1 \le j \le n$, \begin{align*}
    \lambda=\sum_{i=1, i \ne j}^{n}\langle R^{*}(X_i,X_j)X_j,X_i\rangle^{*}&=\sum_{i=1, i \ne j}^{n}K^{*}(X_i,X_j)\\
    &=\sum_{i=1, i \ne j}^{n}[\langle \beta(X_i,X_i), \beta(X_j,X_j)\rangle+\psi(X_i,X_i)+\psi(X_j,X_j)].
\end{align*} Let $f$ be the Example \ref{cylinder}. Substituting the M\"{o}bius invariant tensors of Example \ref{cylinder} into the preceding equation, we obtain \begin{equation}
    \lambda=\lambda\langle X_1,X_1\rangle^{*}=(n-1)\left(-\frac{\kappa_{ss}}{\kappa^{3}}+\frac{\kappa_s^{2}}{k^{4}}\right)\label{eq1}
\end{equation}
and for $2 \le j \le n$
\begin{equation} 
    \lambda=\lambda\langle X_j,X_j\rangle^{*}
=-\frac{\kappa_{ss}}{\kappa^{3}}-(n-3)\frac{\kappa_s^{2}}{\kappa^{4}}.\label{eq2}
\end{equation} From equations \eqref{eq1} and \eqref{eq2}, \begin{equation*}
    \frac{\lambda}{n-1}-\frac{\kappa_s^{2}}{\kappa^{4}}=\lambda+(n-3)\frac{\kappa_s^{2}}{\kappa^{4}}
\end{equation*} which is equivalent to $$\left[\frac{d}{ds}\frac{1}{\kappa}\right]^{2}+\frac{\lambda}{n-1}=0.$$ From the preceding equation, we conclude that $\lambda \leq 0.$ If $\lambda=0$, then $\kappa(s)=\frac{1}{r}$, for all $s\in \mathbb{R}$, where $r>0$ is constant. If $\lambda<0$, then $\kappa(s)=\sqrt{\frac{n-1}{-\lambda}}\frac{1}{s},$ for all $s>0.$

Let $f$ be the Example \ref{cone1}. We have 
\begin{align}
    \lambda=\lambda \langle X_1,X_1\rangle^{*}=(n-1)\left(-\frac{\kappa_{ss}}{\kappa^{3}}+\frac{\kappa_s^{2}}{\kappa^{4}}\right)\label{eq3}
\end{align}
and for $2 \le j \le n,$
\begin{align}
    \lambda=\lambda \langle X_j,X_j\rangle^{*}=-\frac{\kappa_{ss}}{\kappa^{3}}-(n-3)\frac{\kappa_s^{2}}{\kappa^{4}}-\frac{n-2}{\kappa^{2}}. \label{eq4}
\end{align}
It follows from equations \eqref{eq3} and \eqref{eq4} that $$\frac{\lambda}{n-1}-\frac{\kappa_s^{2}}{\kappa^4}=\lambda+(n-3)\frac{\kappa_s^{2}}{\kappa^{4}}+\frac{n-2}{\kappa^{2}}$$ which is equivalent to $$\left[\frac{d}{ds}\frac{1}{\kappa}\right]^{2}+\left[\frac{1}{\kappa}\right]^{2}=-\frac{\lambda}{n-1}.$$ It follows from this equation that $\lambda<0$ and its solution is $\kappa(s)=\sqrt{\frac{n-1}{-\lambda}}\frac{1}{\sin s}$ for all $s \in (0, \pi).$

Let $f$ be the Example \ref{rotacao}. We have 
\begin{align}
    \lambda=\lambda \langle X_1,X_1\rangle^{*}=(n-1)\left(\frac{\kappa_{ss}}{\kappa^{3}}-3\frac{\kappa_s^{2}}{\kappa^{4}}\right)\label{eq5}
\end{align} and for $2 \le j \le n$
\begin{align}
    \lambda=\lambda\langle X_j,X_j\rangle^{*}=\frac{\kappa_{ss}}{\kappa^{3}}-(n+1)\frac{\kappa_s^{2}}{\kappa^{4}}+\frac{(n-2)}{\kappa^{2}}.\label{eq6}
\end{align} It follows from equations \eqref{eq5} and \eqref{eq6} that 
$$\frac{\lambda}{n-1}+3\frac{\kappa_s^{2}}{\kappa^{4}}=\lambda+(n+1)\frac{\kappa_s^{2}}{\kappa^{4}}-\frac{n-2}{\kappa^{2}},$$ which is equivalent to $$-\left[\frac{d}{ds}\frac{1}{\kappa}\right]^{2}+\frac{1}{\kappa^{2}}=\frac{\lambda}{n-1}.$$ The solutions are $$\kappa(s)=\begin{cases}
    \sqrt{\frac{n-1}{\lambda}}\frac{1}{\cosh s}, \quad \lambda>0 \,\,\, s \in \mathbb{R};\\
    \sqrt{\frac{n-1}{-\lambda}}\frac{1}{\sinh s}, \quad \lambda<0\,\,\, s >0;\\
    e^{s}, \quad \lambda=0, \,\,\, s \in \mathbb{R}.
\end{cases} $$

In summary, according to Theorem 4.1 of \cite{AT}, we prove that if $f$ is as Examples \ref{cylinder}, \ref{cone1} and \ref{rotacao} and has constant M\"{o}bius Ricci curvature $\lambda$ then $f$ has constant M\"{o}bius sectional curvature. This proves the statement in item a) of Theorem \ref{main1}.

Now, if $f:M^{n}\to \mathbb{R}^{n+p}$ has flat normal bundle and closed M\"{o}bius form with two distinct principal normal vector fields of multiplicities $k$ and $n-k$, where $2 \le k \le n-2,$ then, according to Lemma \ref{class-two}, $f(M)$ is an open subset of one of the following submanifolds

\begin{description}
    \item[i)] the standard cylinder $\mathbb{S}^{k}\times \mathbb{R}^{n-k}$ in $\mathbb{R}^{n+1};$ 
    \item[ii)] the image of $\sigma^{-1}\circ \tau$ of the hyperbolic cylinder $\mathbb{S}^{k}(r)\times \mathbb{H}^{n-k}(\sqrt{1+r^{2}}) \subset \mathbb{H}^{n+1}$. 
    \item[iii)] the image of $\sigma^{-1}$ of the isoparametric torus $\mathbb{S}^{k}(r)\times \mathbb{S}^{n-k}(\sqrt{1-r^{2}})\subset \mathbb{S}^{n+1}$,
\end{description} where $\sigma:\mathbb{R}^{n+1}\to \mathbb{S}^{n+1}-{(-1,0, \ldots, 0)}$ and $\tau: \mathbb{H}^{n+1}\to \mathbb{S}^{n+1}_{+}$ are conformal diffeomorphism defined by \begin{align*}
\sigma(u)&=\left(\frac{1-||u||^{2}}{1+||u||^{2}}, \frac{2u}{1+||u||^{2}}\right),\\
\tau(y)&=\left(\frac{1}{y_0}, \frac{\bar{y}}{y_0}\right), \quad {y}=(y_0, y_1, \ldots, y_{n+1}):=(y_0, \bar{y}) \in \mathbb{H}^{n+1}.
\end{align*}

If we assume that the examples i), ii), iii) have constant M\"{o}bius Ricci curvature, then it follows from Propositions 3.2, 3.3 and 3.5 of \cite{Ricci} that only the torus $$\mathbb{S}^{k}\left(\sqrt{\frac{k-1}{n-2}}\right)\times \mathbb{S}^{n-k}\left(\sqrt{\frac{n-k-1}{n-2}}\right)$$ for $2 \le k \le n-2$, has constant M\"{o}bius Ricci curvature $$\lambda=\frac{(n-k-1)(n-1)(k-1)}{k(n-k)(n-2)}.$$ This proves the statement in item b) of Theorem \ref{main1}.

\section{The case of at least three distinct principal normals}\label{three}

Let $f:M^{n}\to \mathbb{R}^{m}$, $n \geq 3$, be an umbilical-free proper isometric immersion with flat normal bundle. Let $\eta_1,\ldots, \eta_t$, $t \geq 3,$ be the distinct principal normal vector fields of $f$ with multiplicities $m_1, \ldots, m_t$, respectively, and denote by $E_{\eta_1}, \ldots, E_{\eta_t}$ the corresponding eigenbundle distributions. Assuming that $f$ has constant M\"{o}bius Ricci curvature $\lambda$, then $\omega$ is closed. Therefore, there exists an orthonormal frame $\{X_1, \ldots, X_n\}$ of $(M^{n}, \langle \cdot, \cdot \rangle^{*})$ that diagonalizes $\beta$ and $\psi$ simultaneously, such that for each $1 \le j \le t $ $$\beta(X_i,X_i)=\bar{\eta}_j$$ for all $m_1+\ldots+m_{j-1}+1 \le i \le m_1+\ldots+m_j.$ It follows from Proposition \ref{blaschke2} that for each $1 \le j \le t$, $$(n-2)\psi(X_i,X_i)=\lambda+||\bar{\eta}_j||^{2}-\tr \psi$$ for all $m_1+\ldots+m_{j-1}+1 \le i \le m_1+\ldots+m_j.$ Therefore, denoting by \begin{equation}\label{blaschke-eigenvalue}
\theta_j=\frac{1}{n-2}(\lambda+||\bar{\eta}_j||^{2}-\tr \psi)
\end{equation} for each $1 \le j \le t,$ we see that $\psi$ has at most $t$ eigenvalues $\theta_1, \ldots, \theta_t.$ Since $\lambda$ and $\tr \,\psi$ are constants, we obtain \begin{equation*}\label{fundamental-condition}X(||\bar{\eta}_j||^{2}+(2-n)\theta_j)=0, \quad \forall \,\, X \in \mathfrak{X}(M) \quad \mbox{and} \quad 1 \le j \le t.\end{equation*}

Which are the umbilic-free isometric immersions with flat normal bundle, constant M\"{o}bius Ricci curvature and vanishing M\"{o}bius form that have at least three distinct principal normal vector fields? From now on, we always assume that $\omega$ vanishes identically.

\begin{lemma}\label{dupin}
 If $\nabla^{\perp}_{X}\bar{\eta}_i=0$ for all $X \in E_{\eta_i}$, then $\eta_i$ is a Dupin principal normal vector field.
\end{lemma}

\begin{proof}

    Recall that $$\omega(X)=-\frac{1}{\rho}(\nabla^{\perp}_{X}\mathcal{H}^{f}+\beta(X, \grad^{*}\rho)).$$ Since $\omega=0$ and $\nabla^{\perp}_{X}\bar{\eta}_i=0$ for all $X \in E_{\eta_i}$, we have that \begin{align*}
        0&=\nabla^{\perp}_{X}\mathcal{H}^{f}+\beta(X, \grad^{*}\rho)\\
        &=\nabla^{\perp}_{X}\mathcal{H}^{f}+\rho(\alpha(X, \grad^{*}\rho)-\langle X, \grad^{*}\rho\rangle \mathcal{H}^{f})\\
        &=\nabla^{\perp}_{X}\mathcal{H}^{f}+\rho (\langle X, \grad^{*}\rho \rangle \eta_i-\rho^{-2}X(\rho)\mathcal{H}^{f})\\
        &=\nabla^{\perp}_{X}\mathcal{H}^{f}+\rho^{-1}X(\rho)(\eta_i- \mathcal{H}^{f}).
    \end{align*}
    By definition, $\bar{\eta}_i=\rho^{-1}(\eta_i-\mathcal{H}^{f})$, then \begin{align*}
        0=\nabla^{\perp}_{X}\bar{\eta}_i=\nabla^{\perp}_{X}(\rho^{-1}(\eta_i-\mathcal{H}^{f}))=-\rho^{-2}X(\rho)(\eta_i-\mathcal{H}^{f})+\rho^{-1}(\nabla^{\perp}_{X}\eta_i-\nabla^{\perp}_{X}\mathcal{H}^{f}).
    \end{align*} Therefore, $$\nabla^{\perp}_{X}\eta_i=0$$ for all $X \in E_{\eta_i}.$
\end{proof}

From now on, we write $X \in E_{\eta_i}$ meaning that $X$ is an element of $\{X_1, \ldots, X_n\}$ that generates $E_{\eta_i}$.

\begin{lemma}\label{derivadas}
    Let $X,Y \in E_{\eta_i}$ and $Z \in E_{\eta_j}$ with $i \ne j$, then 
    \begin{description}
    \item[a)] $\langle X,Y\rangle^{*}\nabla^{\perp}_{Z}\bar{\eta}_i=\langle \nabla^{*}_{X}Y,Z\rangle^{*}(\bar{\eta}_i-\bar{\eta}_j).$
    \item[b)] $\langle X,Y\rangle^{*}Z(\theta_i)=\langle \nabla^{*}_{X}Y,Z\rangle^{*}(\theta_i-\theta_j).$ 
    \item[c)] $\langle X,Y\rangle^{*}\frac{Z(||\bar{\eta}_i||^{2}+2\theta_i)}{2}=\langle \nabla^{*}_{X}Y,Z\rangle^{*}(||\bar{\eta}_i||^{2}+2\theta_i-K^{*}(X,Z)).$ 
    \end{description}
\end{lemma}

\begin{proof}
    \begin{description}
        \item[a)]Use the conformal Codazzi equation $$(\nabla^{\perp}_{X}\beta)(Z,Y)=(\nabla^{\perp}_{Z}\beta)(X,Y).$$
        \item[b)] Use the conformal Codazzi equation $$(\nabla^{*}_{X}\psi)(Z,Y)=(\nabla^{*}_{Z}\psi)(X,Y).$$ 
        \item[c)] Taking the inner product of a) with $\bar{\eta}_i$ and summing it with b) we get $$\langle X,Y\rangle^{*}\frac{Z(||\bar{\eta}_i||^{2}+2\theta_i)}{2}=\langle \nabla^{*}_{X}Y, Z\rangle^{*}(||\bar{\eta}_i||^{2}+\theta_i-\langle \bar{\eta}_i, \bar{\eta}_j\rangle-\theta_j).$$ Using the equation $$K^{*}(X,Z)=\langle \bar{\eta}_i, \bar{\eta}_j\rangle+\theta_i+\theta_j$$ we obtain c).\end{description}\end{proof}

\begin{lemma}\label{totally-geodesic}
    The eigenbundle distribution $E_{\eta_i}$ is totally geodesic for all $1 \le i \le t.$
\end{lemma}

\begin{proof}
    Substituting \eqref{blaschke-eigenvalue} in item c) of Lemma \ref{derivadas} we get $$\frac{n}{2}\langle X,Y\rangle^{*}Z(\theta_i)=\langle \nabla^{*}_{X}Y,Z\rangle^{*}(n \theta_i-\lambda+\tr \, \psi-K^{*}(X,Z)).$$ However, substituting item b) of Lemma \ref{derivadas} into the preceding equation, we obtain $$\langle \nabla^{*}_{X}Y,Z\rangle^{*}\left(\frac{n}{2}(\theta_i+\theta_j)-\lambda+\tr \psi-K^{*}(X,Z)\right)=0.$$ If, for some $x \in M^{n}$, $K^{*}(X,Z)=\frac{n}{2}(\theta_i+\theta_j)-\lambda+\tr \, \psi$, then using $K^{*}(X,Z)=\langle \bar{\eta}_i, \bar{\eta}_j\rangle+\theta_i+\theta_j$ we obtain $$\langle \bar{\eta}_i, \bar{\eta}_j \rangle+\frac{2-n}{2}(\theta_i+\theta_j)=-\lambda+\tr\, \psi.$$ Substituting \eqref{blaschke-eigenvalue} into the preceding equation, we get $$\langle \bar{\eta}_i, \bar{\eta}_j\rangle-\frac{2\lambda+||\bar{\eta}_i||^{2}-2\tr\, \psi+||\bar{\eta}_j||^{2}}{2}=-\lambda+\tr\, \psi$$ which is equivalent to $$||\bar{\eta}_i-\bar{\eta}_j||=0.$$ Therefore, $\bar{\eta}_i=\bar{\eta}_j,$ which is a contradiction. This proves that $$\langle \nabla^{*}_{X}Y,Z\rangle^{*}=0$$ on $M^{n}$ for all $X,Y \in E_{\eta_i}$ and all $Z \in E_{\eta_j}$ with $i \ne j.$
\end{proof}

We recall that an umbilic-free isometric immersion $f:M^{n}\to \mathbb{R}^{m}$ is said to be M\"{o}bius isoparametric if its M\"{o}bius form vanishes identically and its M\"{o}bius principal normal vector fields are parallel with respect to the normal connection of $f$ (see \cite{A}).

\begin{lemma}\label{isoparametric}
    If $f:M^{n}\to \mathbb{R}^{m}$, $n \geq 3$, is an isometric immersion with constant M\"{o}bius Ricci curvature and flat normal bundle with vanishing M\"{o}bius form then $f$ is M\"{o}bius isoparametric.
\end{lemma}

\begin{proof}
  From Lemmas \ref{derivadas}-a) and \ref{totally-geodesic}  follow that $\nabla^{\perp}_{Z}\bar{\eta}_i=0$ for all $Z \in E_{\eta_j}$  and all $1 \le j \ne i \le t.$ 
    Since $\beta$ is traceless, we have $\sum_{i=1}^{t}m_i \bar{\eta}_i=0.$ Then, for all $Z \in E_{\eta_i}$, we get $\nabla^{\perp}_{Z}\bar{\eta}_i=0.$  Therefore, $\bar{\eta}_i$, for all $1 \le i \le t$, is parallel with respect to the normal connection of $f$.
\end{proof}

\begin{lemma}\label{holonomic}
At any point $x\in M^{n}$, $$\langle \nabla^{*}_{X}Y,Z\rangle^{*}=0$$ for all $X \in E_{\eta_i}$, $Y\in E_{\eta_j}$ and $Z \in E_{\eta_\ell},$ with $i \ne j \ne \ell \ne i.$
\end{lemma}

\begin{proof}
   Let $F=\rho(\Psi \circ f)$ be the M\"{o}bius lift of $f$. The second fundamental form of $F$ is given by $$\alpha^{F}(X,Y)=\phi(\beta(X,Y))-\psi(X,Y)F-\langle X,Y\rangle^{*}\zeta$$ for all $X,Y \in \mathfrak{X}(M)$ (see \eqref{sf}). A computation of the Ricci tensor of $(M^{n}, \langle \cdot, \cdot \rangle^{*})$ using the Gauss equation of $F$ yields \begin{equation}\label{ricci-tensor}
   \Ric^{*}(X,Y)=n\langle \alpha^{F}(X,Y), \mathcal{H}^{F}\rangle-\sum_{i=1}^{n}\langle \beta(X,X_i),\beta(Y,X_i)\rangle-2\psi(X,Y)\end{equation} for all $X,Y \in \mathfrak{X}(M).$

   For all $X \in E_{\eta_i}$, define $\eta_i^{*}$ by $$\eta_i^{*}=\phi(\bar{\eta}_i)-\theta_iF-\zeta.$$ Therefore, the principal normal vector fields of $F$ are $\eta_i^{*}$ for $1 \le i \le t$ whose associated eigenbundle distributions are $E_{\eta_i}.$ It is clear that $f$ has constant M\"{o}bius Ricci curvature $\lambda$ if and only if $F$ has Ricci curvature $\lambda.$ Using \eqref{ricci-tensor}, we get that $$-\lambda=\langle \eta_i^{*}, \eta_i^{*}\rangle-n\langle \mathcal{H}^{F}, \eta_i^{*}\rangle, \quad 1 \le i \le t.$$ This in turn is equivalent to $$\hat{\eta}_i=\eta_i^{*}-\frac{n}{2}\mathcal{H}^{F}$$  satisfy \begin{equation}\label{norm-hat}
   \langle \hat{\eta}_i, \hat{\eta}_i \rangle=\frac{n^{2}}{4}\langle \mathcal{H}^{F}, \mathcal{H}^{F}\rangle-\lambda\end{equation} for all $1 \le i \le t.$

   We claim that at any point the vectors $\eta_i^{*}-\eta_j^{*}$ and $\eta_i^{*}-\eta_{\ell}^{*}$ are linearly independent if $i \ne j \ne \ell \ne i.$ Assume to the contrary that $$\eta_i^{*}-\eta_j^{*}=\mu(\eta_i^{*}-\eta_\ell^{*})$$ for some $\mu \ne 0.$ Then $$(1-\mu)\hat{\eta}_i=\hat{\eta}_j-\mu \hat{\eta}_\ell$$ yields $$\langle \hat{\eta}_i, \hat{\eta}_i \rangle-2\mu \langle \hat{\eta}_i, \hat{\eta}_i\rangle+ \mu^{2}\langle \hat{\eta}_i, \hat{\eta}_i\rangle=\langle \hat{\eta}_j, \hat{\eta}_j \rangle -2\mu \langle \hat{\eta}_j, \hat{\eta}_\ell \rangle+\mu^{2}\langle \hat{\eta}_\ell, \hat{\eta}_\ell \rangle.$$ Using \eqref{norm-hat} we obtain $$\langle \hat{\eta}_j, \hat{\eta}_\ell \rangle=\langle \hat{\eta}_i, \hat{\eta}_i \rangle=\langle \hat{\eta}_j, \hat{\eta}_j \rangle=\langle \hat{\eta}_\ell, \hat{\eta}_\ell \rangle.$$ Then $$\langle \hat{\eta}_j, \hat{\eta}_j-\hat{\eta}_{\ell} \rangle=0 \quad \mbox{and} \quad \langle \hat{\eta}_{\ell}, \hat{\eta}_\ell- \hat{\eta}_j \rangle=0,$$ and hence $\langle \hat{\eta}_j-\hat{\eta}_\ell, \hat{\eta}_j-\hat{\eta}_\ell \rangle=0.$ However, this is equivalent to $$\langle \phi(\bar{\eta}_j)-\theta_jF-\phi(\bar{\eta}_\ell)+\theta_\ell F, \phi(\bar{\eta}_j)-\theta_j F-\phi(\bar{\eta}_\ell)+\theta_\ell F \rangle=0.$$ Therefore, using that $\phi$ is an isometry and $F$ is a light-like vector, we get $$|| \bar{\eta}_j-\bar{\eta}_\ell||=0,$$ that is $\bar{\eta}_j=\bar{\eta}_\ell,$ which is a contradiction. Thus, $\eta_i^{*}-\eta_j^{*}$ and $\eta_i^{*}-\eta_\ell^{*}$ are linearly independent vectors.

   By the Codazzi equation of $F$, we get $$\langle \nabla^{*}_{X}Y,Z\rangle^{*}(\eta_\ell^{*}-\eta_i^{*})=\langle \nabla^{*}_{Y}X,Z\rangle^{*}(\eta_j^{*}-\eta_i^{*}),$$ then $\langle \nabla^{*}_{X}Y,Z\rangle^{*}=\langle \nabla^{*}_{Y}X,Z\rangle^{*}=0$ for all $X \in E_{\eta_j}$, $Y \in E_{\eta_\ell}$ and $Z \in E_{\eta_i}$ with $1 \le i \ne j \ne \ell \ne i \le t$.
    \end{proof}

\begin{lemma}\label{totally-geodesic2}
    The smooth distribution $E_{\eta_i}^{\perp}$ is totally geodesic for all $1 \le i \le t.$
\end{lemma}

\begin{proof}
    This follows immediately from Lemmas \ref{totally-geodesic} and \ref{holonomic}.
\end{proof}

In the following, we build some standard examples of isometric immersions with flat normal bundle and at least three distinct principal normal vector fields.

\begin{example}\label{cilindro} 
						Let $g:M^{n-p}\to \mathbb{R}^{m-p}$, $1 \le p \le n-2,$ be an umbilic-free isometric immersion with $R^{\perp} \equiv 0$.  The cylinder in $\mathbb{R}^{m}$ over $g$ is the immersion given by \begin{equation*}
							f=g \times \Id:M^{n-p}\times \mathbb{R}^{p}\to \mathbb{R}^{m-p}\times \mathbb{R}^{p}=\mathbb{R}^{m},
						\end{equation*} where $\Id:\mathbb{R}^{p}\to \mathbb{R}^{p}$ is the identity map.

						Since the induced metric by $f$ is \begin{equation*}\label{metric1}\langle \,\,,\,\rangle_{f}=ds^{2}+du_{n-p+1}^{2}+\ldots+du_{n}^{2},
						\end{equation*}
						where $(u_{n-p+1}, \ldots, u_n)$ are the canonical coordinates on $\mathbb{R}^{p}$, the M\"{o}bius metric determined by $f$ is  \begin{equation*}\label{M\"{o}biuseuclideo}\langle\,\,,\,\rangle^{*}=\rho^{2}(ds^{2}+du_{n-p+1}^{2}+\ldots+du_n^{2}).\end{equation*}
					\end{example}

                    \begin{example}\label{coned}
						Let $g:M^{n-p}\to \mathbb{S}^{m-p}$, $1 \le p \le n-2$, be an umbilic-free isometric immersion with $R^{\perp}\equiv 0$.
						The generalized cone in $\mathbb{R}^{m}$ over $g$ is the immersion $f:M^{n-p}\times \mathbb{H}^{p}\to \mathbb{R}^{m}$ given by \begin{equation*} f(x,z)=\Theta \circ (g,\Id)(x,z)=(z_1g(x), z_2, \ldots, z_{p}),\end{equation*} where $\mathbb{H}^{p}=\mathbb{R}_{+}\times \mathbb{R}^{p-1}$ is endowed with the hyperbolic metric $$dz^{2}=\frac{1}{z_1^{2}}(dz_1^{2}+\ldots+dz_{p}^{2}),$$ $\Id:\mathbb{H}^{p}\to \mathbb{H}^{p}$ is the identity map and $\Theta:\mathbb{S}^{m-p}\times \mathbb{H}^{p}\to \mathbb{R}^{m}$ is the conformal diffeomorphism defined by $$ \Theta(y,z)=(z_1y, z_2, \ldots, z_{p}),$$ whose conformal factor is $z_1.$ Then, the induced metric by  $f$ is \begin{equation*}
							\langle \,\,,\,\rangle_{f}=z_1^{2}(ds^{2}+dz^{2}).
						\end{equation*}

						 The second fundamental form of $f$ is given by $$\alpha^{f}\left(\frac{X_i}{z_1}, \frac{X_j}{z_1}\right)=\frac{1}{z_1} \alpha^{g}(X_i,X_j), \quad \alpha^{f}\left(\frac{d}{dz_1}, \frac{d}{dz_1}\right)=0=\alpha^{f}\left(\frac{d}{dz_1}, \frac{X_i}{z_1}\right),$$ where $X_1, \ldots, X_{n-p}$ is an orthonormal frame with respect to $ds^{2}$ that diagonalizes $\alpha^{g}.$  Thus the M\"{o}bius metric determined by $f$ is \begin{equation*}\label{M\"{o}bius-metric-esfera}
							\langle\,\,,\,\rangle^{*}=(\rho z_1)^{2}(ds^{2}+dz^{2}).
						\end{equation*}

						\end{example}

\begin{example}\label{subvariedade-de-rotacao}
							Let $g:M^{n-p}\to \mathbb{H}^{m-p}$, $1 \le p \le n-2$, be an umbilic-free isometric immersion with $R^{\perp}\equiv 0$ where $$\mathbb{H}^{m-p}=\{(z_1, \ldots, z_{m-p}): z_{m-p}>0\} $$ is endowed with the hyperbolic metric $$ ds^{2}=\frac{1}{z_{m-p}^{2}}(dz_1^{2}+\ldots+dz_{m-p}^{2}).$$

							The rotational submanifold over $g$ is the immersion $f:M^{n-p}\times \mathbb{S}^{p}\to \mathbb{R}^{m}$ defined by \begin{equation*}\label{rotational} f(x,y)=\Theta \circ (g, \Id)(x,y)=(z_1(x),\ldots,z_{m-p-1}(x),z_{m-p}(x)y),\end{equation*} where $\Id:\mathbb{S}^{p}\to \mathbb{S}^{p}$ is the identity map and $\Theta: \mathbb{H}^{m-p}\times \mathbb{S}^{p}\to \mathbb{R}^{m}$ is the conformal diffeomorphism $$\Theta(z,y)=(z_1, \ldots, z_{m-p-1}, z_{m-p}y),$$ whose conformal factor is $z_{m-p}$. Then, the induced metric by $f$ is $$ \langle \,\,,\,\rangle_{f}=z_{m-p}^{2}(x)(g^{*}ds^{2}+dy^{2}|_{\mathbb{S}^{p}}),$$ $x \in M^{n-p}$.

				 	The second fundamental form of $f$ is given by \begin{align*}
								\alpha^{f}\left(\frac{X_i}{z_{m-p}}, \frac{X_j}{z_{m-p}}\right)&=\begin{cases}
									\frac{1}{z_{m-p}^{2}}\Theta_{*}\alpha^{g}({X_i}, {X_i})-\frac{1}{z_{m-p}^{3}}\Theta_{*}((\grad \, z_{m-p}) \circ g)^{\perp}, \quad i=j,\\
									0, \quad i \ne j,
								\end{cases}\\
								\alpha^{f}\left(\frac{1}{z_{m-p}}e_k, \frac{1}{z_{m-p}}e_r\right)&=\begin{cases}-\frac{1}{z_{m-p}^{3}}\Theta_{*}((\grad\,z_{m-p})\circ g)^{\perp}, \quad k=r,\\
									0, \quad k \ne r,
								\end{cases}
							\end{align*} where $X_1, \ldots, X_{n-p}$ and $e_{n-p+1}, \ldots, e_{n}$ are orthonormal frames with respect to $g^{*}ds^{2}$ and $dy|_{\mathbb{S}^{p}}$, respectively. Here $\grad$ is computed with respect to $ds^{2}$. The description of the second fundamental form of $f$ can be found in Exercise 9.2 of \cite{DT}.

									The M\"{o}bius metric determined by $f$ is given by  $$ \langle\,\,,\,\rangle^{*}=(\rho z_{m-p}(x))^{2}(g^{*}ds^{2}+dy^{2}|_{\mathbb{S}^{p}}),$$ $x\in M^{n-p}.$
								\end{example} 
					The following observation will be useful for the proof of Theorem \ref{main2}. Let $f:M \to \tilde{M}$ be an immersion and $g_1$ and $g_2$ are metrics on $\tilde{M}$ such that $g_2=\lambda^{2}g_1$ for some positive function $\lambda \in C^{\infty}(\tilde{M})$, then $$f_j=f:(M, f^{*}g_j)\to (\Tilde{M},g_j),$$ $1 \le j \le 2$, are isometric immersions whose induced metrics satisfy $f^{*}g_2=(\lambda \circ f)^{2}f^{*}g_1.$ Denote $^{f_1}\tilde{\nabla}$ and $^{f_2}\Tilde{\nabla}$ the connections of the induced vetor bundle $f_i^{*}T\tilde{M}$, $1 \le i \le 2$, then taking the normal component of $${^{f_2}\tilde{\nabla}_{X}}\xi={^{f_1}\tilde{\nabla}_{X}}\xi+\frac{1}{\lambda}(X(\lambda)\xi+\xi(\lambda)X),$$ for $X\in f_{*}TM$ and $\xi \in N_{f}M$, we obtain \begin{equation}\label{con-conf}
									{^{f_2}\nabla^{\perp}_{X}}\xi={^{f_1}\nabla^{\perp}_{X}}\xi+\frac{1}{\lambda}X(\lambda)\xi.
								\end{equation}

\begin{theorem}\label{main2}
    Let $f:M^{n}\to \mathbb{R}^{m}$, $n \geq 3$, be an isometric immersion with constant M\"{o}bius Ricci curvature, flat normal bundle and vanishing M\"{o}bius form. If $f$ has at least three distinct principal normal vector fields, then $f$ has vanishing M\"{o}bius sectional curvature and hence it is M\"{o}bius equivalent to one of the following submanifolds:
    
\begin{description}
\item[a)]$\mathbb{S}^{1}(r_1)\times \ldots \times \mathbb{S}^{1}(r_{n-s})\times \mathbb{R}^{s}\subset \mathbb{R}^{2n-s},$ for all $1 \le s \le n-2$.
\item[b)] the three dimensional cone over the flat torus $\mathbb{S}^{1}(r)\times \mathbb{S}^{1}(\sqrt{1-r^{2}})\subset \mathbb{S}^{3}.$
\end{description}
\end{theorem}

\begin{proof}
    Let $\eta=\eta_t$ be a principal normal vector field of multiplicity $p \geq 1.$ By Lemmas \ref{dupin} and \ref{isoparametric}, $\eta$ is a Dupin vector field.

   By Lemma \ref{totally-geodesic2} we see that $E_{\eta}^{\perp}$ is totally geodesic with respect to $\langle \cdot, \cdot \rangle^{*}$, then from Proposition 4.1 of \cite{AT} we see that $E_{\eta}^{\perp}$ is spherical with respect to the metric induced by $f$. The main theorem of \cite{DFT} implies that $f$ is locally, up to a composition with a conformal transformation in $\mathbb{R}^{m}$, an isometric immersion $f=\Theta \circ (g, \Id):M^{n-p}\times \mathbb{Q}^{p}_{-c}\to \mathbb{R}^{m}$, where $g:M^{n-p}\to \mathbb{Q}^{m-p}_{c}$, $1 \le p  <n-1,$ is an umbilic-free isometric immersion with flat normal bundle and $\Theta=\Id:\mathbb{R}^{m}\to \mathbb{R}^{m}$ if $c=0$ or $\Theta$ is given in Examples \ref{coned} and \ref{rotational} if $c =1$ or $c=-1$, respectively.

   Since we have shown in Examples \ref{cilindro}-\ref{coned} and \ref{subvariedade-de-rotacao}, $\eta=0$ or $\eta=-\frac{1}{z_{m-p}^{3}}\Theta_{*}((\grad\,z_{m-p})\circ g)^{\perp}$, respectively.

   It follows from Lemmas \ref{totally-geodesic}, \ref{totally-geodesic2} that the distributions $E_{\eta}$ and $E_{\eta}^{\perp}$ are totally geodesic with respect to $\langle\,\,,\,\rangle^{*}.$ In particular, due to the de Rham's theorem, $(M^{n}, \langle \, \,,\,\rangle^{*})$ is locally a Riemannian product. 
									Since the M\"{o}bius metric determined by $f$ has the form \begin{equation}\label{product-metric}\langle \,\,,\,\rangle^{*}=\theta^{2}(ds^{2}+\sigma_{-c}),\end{equation} where $ds^{2}$ is the induced metric by $g$, $\sigma_{-c}$ is the canonical metric on $\mathbb{Q}^{p}_{-c}$ and $$\theta=\begin{cases}
									    \rho, \quad \mbox{if}\,\,\, c=0,\\
                                        \rho z_1, \quad \mbox{if} \,\,\, c=1,\\
                                        \rho z_{m-p}(x), \,\, x \in M^{n-p}, \quad \mbox{if} \,\, c=-1,
									\end{cases}$$ then $\theta$ is a non-vanishing constant function.

                                    Now, we prove that the mean curvature vector field of $g:U \subset M^{n-p} \to \mathbb{Q}^{m-p}_{c}$ is parallel with respect to $^{g}\nabla^{\perp}$, for any $c=0,1$ or $-1.$ Observe that since $f$ is M\"{o}bius isoparametric, the M\"{o}bius principal normal vector field $\bar{\eta}$ is parallel along the normal connection of $f.$

                                    Let $f$ be the cylinder over $g$. Then its M\"{o}bius principal normal vector field associated with $\eta=0$ is $$\bar{\eta}=-\rho^{-1}\mathcal{H}^{f}=-\frac{n-p}{n}\rho^{-1}\mathcal{H}^{g}.$$ Hence $^{g}\nabla^{\perp}\mathcal{H}^{g} \equiv 0.$

                                    Let $f$ be the generalized cone over $g$.  Then its M\"{o}bius principal normal vector field associated with $\eta=0$ is $$ \bar{\eta}=-\rho^{-1}\mathcal{H}^{f}=-\frac{n-p}{n}\frac{1}{\rho z_1}\mathcal{H}^{g}.$$
									Therefore, ${^{g}\nabla^{\perp} \mathcal{H}^{g} }\equiv 0.$

                                    Let $f$ be the rotational submanifold over $g$. Since $\eta=-\frac{1}{z_{m-p}^{3}}\Theta_{*}((\grad z_{m-p})\circ g)^{\perp}$ we have  $\eta-\mathcal{H}^{f}=-\frac{1}{n}\frac{1}{z_{m-p}^{2}}(n-p)\Theta_{*}\mathcal{H}^{g}, $ then \begin{equation}
										\bar{\eta}=\rho^{-1}(\eta-\mathcal{H}^{f})=-\frac{n-p}{n}\frac{1}{\rho z_{m-p}}\Theta_{*}\frac{\mathcal{H}^{g}}{z_{m-p}}.
									\end{equation}
									Using \eqref{con-conf} we obtain \begin{align*}
										0={^{f}\nabla^{\perp}_{X}}\Theta_{*}\frac{\mathcal{H}^{g}}{z_{m-p}}&=\Theta_{*}{^{g}\nabla^{\perp}_{X}}\frac{\mathcal{H}^{g}}{z_{m-p}}+\frac{X(z_{m-p})}{z_{m-p}^{2}}\Theta_{*}\mathcal{H}^{g}\\
										&=\frac{1}{z_{m-p}}\Theta_{*}{^{g}\nabla^{\perp}_{X}}\mathcal{H}^{g},
									\end{align*} for every $X \in \mathfrak{X}(M^{n-p}).$ Therefore, $^{g}\nabla^{\perp}\mathcal{H}^{g}\equiv 0.$

                                  By hypothesis we know that the product metric \eqref{product-metric} has constant Ricci curvature $\lambda.$ Then it follows from the description of the curvature tensor of the product metric that the isometric immersion $g:U\subset M^{n-p}\to \mathbb{Q}^{m-p}_{c}$ has constant Ricci curvature.

                                  In summary, we prove that $f(M)$ is locally, M\"{o}bius equivalent, to a cylinder, a generalized cone or a rotational submanifold over any isometric immersion $g:U \subset M^{n-p}\to \mathbb{Q}^{m-p}_{c}$ with flat normal bundle, parallel mean curvature vector field and constant Ricci curvature, according $c=0, 1$ or $-1$, respectively.

                                  The conclusion of this theorem for $n-p \geq 3$ follows from Theorem 1 of \cite{CR} and the fact that the M\"{o}bius metric is a product Riemannian metric.

If $n-p=2$ then $g:M^{2} \to \mathbb{Q}^{m-p}_{c}$ is a surface with flat normal bundle. Since $g$ has constant Ricci curvature then $M^{2}$ has constant sectional curvature $\tilde{c}$. Due to the fact that $g:M^{2}_{\tilde{c}}\to \mathbb{Q}^{m-p}_{c}$ has parallel mean curvature vector field, from Theorem 1 of \cite{DT1} we have that one of the following possibilities holds:

\begin{description}
\item[i)] $\tilde{c}=0=c$ and $g(M^{2})$ is an open subset of $$ \mathbb{S}^{1}(r_1)\times \mathbb{S}^{1}(r_2)\subset \mathbb{R}^{4} \quad \mbox{or} \quad \mathbb{S}^{1}(r_1)\times \mathbb{R}\subset \mathbb{R}^{3}.$$
\item[ii)] $c=1, \tilde{c}=0$ and $g(M^{2})$ is an open subset of the torus $$ \mathbb{S}^{1}(r)\times \mathbb{S}^{1}(\sqrt{1-r^{2}}) \subset \mathbb{S}^{3}.$$
\item[iii)] $c=-1, \tilde{c}=0$ and $g(M^{2})$ is an open subset of the hyperbolic cylinder $$\mathbb{H}^{1}(\sqrt{1+r^{2}})\times \mathbb{S}^{1}(r)\subset \mathbb{H}^{3}.$$
\end{description}

Assume $p=1$. If $g:M^{2}\to \mathbb{R}^{3}$ is a round cylinder $\mathbb{S}^{1}(r_1)\times \mathbb{R}$, then the cylinder $f=g \times \Id:M^{2}\times \mathbb{R}\to \mathbb{R}^{4}$ has only two distinct principal curvatures, which is in contradiction to our assumption.

If $f$ is a rotational hypersurface over a hyperbolic cylinder $\mathbb{H}^{1}(\sqrt{1+r^{2}})\times \mathbb{S}^{1}(r)\subset \mathbb{H}^{3}$, then $f:\mathbb{H}^{1}(\sqrt{1+r^{2}})\times \mathbb{S}^{1}(r)\times \mathbb{S}^{1}\to \mathbb{R}^{4}$ is defined by $$f(y_1,y_2,x_1,x_2)=\left(\frac{x_1}{y_1+y_2}, \frac{x_2}{y_1+y_2}\right),$$ where $-y_1^{2}+y_2^{2}=-1-r^{2}$, $y_1>0$, and $x_1:\mathbb{S}^{1}(r)\to \mathbb{R}^{2}$, $x_2:\mathbb{S}^{1}\to \mathbb{R}^{2}$ are umbilical inclusions. Therefore, a rotational hypersurface over a hyperbolic cylinder is M\"{o}bius equivalent to a cone over a homogeneous torus $\mathbb{S}^{1}(r)\times \mathbb{S}^{1}(\sqrt{1-r^{2}})\subset \mathbb{S}^{3}$.

Now, assume $p>1$. Recall that a Riemannian product of constant sectional curvature manifolds $M^{p}_{c_1}\times M^{n-p}_{c_2}$ is an Einstein manifold if and only if $$(p-1)c_1=(n-p-1)c_2.$$ Since $(M^{2}_{0} \times \mathbb{Q}^{p}_{-c}, ds^{2}+\sigma_{-c})$ is an Einstein manifold, then $c=0$. Therefore, $g(M^{2})$ is an open subset of $\mathbb{S}^{1}(r_1)\times \mathbb{S}^{1}(r_2)\subset \mathbb{R}^{4}$.

In summary, if $f:M^{n}\to \mathbb{R}^{m}$, $n-p=2$, is an umbilic-free immersion with constant M\"{o}bius Ricci curvature, flat normal bundle and vanishing M\"{o}bius form with three distinct principal normals then $f(M^{n})$ is locally M\"{o}bius equivalent to a cylinder $$\mathbb{S}^{1}(r_1)\times \mathbb{S}^{1}(r_2)\times \mathbb{R}^{p} \subset \mathbb{R}^{p+4}$$ or a three dimensional cone over the standard torus $\mathbb{S}^{1}(r)\times \mathbb{S}^{1}(\sqrt{1-r^{2}})\subset \mathbb{S}^{3}$. Thus, we complete the proof of the theorem.
\end{proof}

\section*{Acknowledgements}

Mateus Antas is supported by CAPES, grant 88887.286908/2026-00.

\Addresses

\end{document}